\documentclass{birkjour}%
\usepackage{amsmath}
\usepackage{amsfonts}
\usepackage{amssymb}
\usepackage{graphicx}%
\providecommand{\U}[1]{\protect \rule{.1in}{.1in}}
\newtheorem{theorem}{Theorem}

\begin{document}

\title[Sharp $L^2$ decay rate for degenerate OIO's]{\vspace*{-3.0cm} 	\\
	Sharp $L^2$ decay rate for (1+2)-dimensional  oscillatory integral operators with 	cubic  polynomial phases}
\author{Jayden Lang}\address{American Heritage School\\Plantation, FL 33325\\ USA}\email{jaydenlang88@gmail.com}
\author{Wan Tang}\address{Department of Biostatistics and Data Science\\Tulane University\\New Orleans, LA 70112\\USA}\email{wtang1@tulane.edu}
\maketitle

\begin{abstract}
	In this paper, we consider the (1+2)-dimensional oscillatory integral with
	degenerate cubic homogeneous polynomial phase. We prove that the $L^{2}$ decay
	rate of 3/8 given in (Archiv der Mathematik, 122: 437-447, 2024) is sharp.
	
\end{abstract}

\section{Introduction}

Consider the $(1+2)$-dimensional oscillatory integral operators with cubic
homogeneous polynomial phases
\begin{equation}
T_{\lambda}f(u,v)=\int_{\mathbb{R}}e^{i\lambda \left[  P_{1}(u,v)t^{2}%
+P_{2}(u,v)t\right]  }\Phi(u,v,t)f(t)\mathrm{d}t \label{a}%
\end{equation}
where $\lambda>0$ and $P_{1}(u,v)$ ($P_{2}(u,v)$) is a non-zero homogeneous
polynomial of degree $1$ ($2)$. If the phase function $S(u,v,t)=P_{1}%
(u,v)t^{2}+P_{2}(u,v)t$ is non-degenerate in the sense of \cite{tang06}, i.e.,
$P_{2}(u,v)$ has no multiple factors when it is factored into linear terms
over $\mathbb{C}$ and $P_{1}(u,v)$ is non zero, then one has the following
nearly sharp (up to a power of $\log(\lambda)$) estimate:
\[
\left \Vert T_{\lambda}f\right \Vert _{L^{2}\left(  \mathbb{R}^{2}\right)  }\leq
C_{\Phi}\lambda^{-\frac{1}{2}}\log \lambda \Vert f\Vert_{L^{2}(\mathbb{R})},
\]
where $C_{\Phi}$ is a constant not depending on $\lambda$ \cite{tang06}$.$ If
$P_{2}(u,v)$ has multiple factors, then by an affine transformation, the
oscillatory integral operators (\ref{a}) can be transformed to the format
\begin{equation}
T_{\lambda}f(u,v)=\int_{\mathbb{R}}e^{i\lambda \left[  ut^{2}+au^{2}t\right]
}\Phi(u,v,t)f(t)\mathrm{d}t \label{b}%
\end{equation}
if $P_{1}(u,v)$ is a factor of $P_{2}(u,v)$, or%
\begin{equation}
T_{\lambda}f(u,v)=\int_{\mathbb{R}}e^{i\lambda \left[  ut^{2}+v^{2}t\right]
}\Phi(u,v,t)f(t)\mathrm{d}t \label{c}%
\end{equation}
if $P_{1}(u,v)$ is not a factor of $P_{2}(u,v).$

The oscillatory integral operators (\ref{b}) is essentially a (1+1)
dimensional\ case, which has been well studied. The oscillatory integral
operators (\ref{c}) was studied in \cite{tan24}, where the following decay
estimate was established:
\begin{equation}
\left \Vert T_{\lambda}f\right \Vert _{L^{2}\left(  \mathbb{R}^{2}\right)  }\leq
C_{\Phi}\lambda^{-3/8}\Vert f\Vert_{L^{2}(\mathbb{R})}. \label{d}%
\end{equation}
However, \cite{tan24} did not determine the sharpness of this bound, showing
only that the decay rate cannot exceed $\lambda^{-1/2}$. In this paper, we
establish that the decay rate $\lambda^{-3/8}$ is indeed sharp for the
oscillatory integral operator (\ref{c}).

\section{Sharpness of the decay rate}

In this section, we prove the following theorem which, together with
(\ref{d}), establishes the sharpness of the decay rate $\lambda^{-3/8}$ for
the oscillatory integral operator (\ref{c}).

\begin{theorem}
[Sharpness of the decay rate]Let $T_{\lambda}$ be the oscillatory integral
operator defined in \eqref{c}. Then for any amplitude functions $\Phi$ with
$\Phi \left(  0,0,0\right)  \neq0$, we have%
\begin{equation}
\left \Vert T_{\lambda}\right \Vert _{L^{2}\left(  \mathbb{R}\right)
\rightarrow L^{2}\left(  \mathbb{R}^{2}\right)  }\geq C_{\Phi}\lambda^{-3/8}.
\label{t}%
\end{equation}

\end{theorem}

In the following, we use the notation $a\lesssim b$ $\left(  \text{resp.
}a\gtrsim b\right)  $ to mean that $|a|\leq C|b|$ (resp. $|a|\geq C|b|$) for
some constant $C>0$ independent of $\lambda$.

\begin{proof} Since $\Phi(0,0,0)\neq0$, we may assume without loss of
generality that
\[
\Phi(u,v,t)\geq1\quad \text{on}\quad \{(u,v,t):|u|<2r,\ |v|<2r,\ |t|<2r\}
\]
for some $r>0$. If $r<1$, we rescale the variables by setting
\[
u^{\prime}=\frac{u}{r},\qquad v^{\prime}=\frac{v}{r},\qquad t^{\prime}%
=\frac{t}{r}.
\]
Under this change of variables, the amplitude function satisfies
\[
\Phi(u,v,t)\geq1\quad \text{on}\quad \{(u,v,t):|u|<2,\ |v|<2,\ |t|<2\}.
\]
This rescaling introduces a constant factor $D=r^{3}$ into the phase function,
so that
\[
S(u,v,t)=D\bigl(ut^{2}+v^{2}t\bigr),
\]
where $D\leq1$ is a constant.

By the continuity of the phase function $S\left(  u,v,t\right)  =D\left(
ut^{2}+v^{2}t\right)  $, for a small positive number $\delta \left(  <D\right)
$, there is an $\varepsilon>0$ such that for any $\left(  u,v,t\right)  $ with
$u\in(1-\varepsilon,1+\varepsilon),v\in(1-\varepsilon,1+\varepsilon)$ and
$t\in(1-\varepsilon,1+\varepsilon)$, we have $\left \vert S\left(
u,v,t\right)  -S(1,1,1)\right \vert =\left \vert S\left(  u,v,t\right)
-2D\right \vert <\delta$. Choose a smooth function $0\leq f_{\lambda}(t)\leq1$
such that
\[
f_{\lambda}(t)=%
\begin{cases}
1, & \lambda^{1/2}t\in \left(  1-\varepsilon/2,1+\varepsilon/2\right) \\
0, & \lambda^{1/2}t\notin \left(  1-\varepsilon,1+\varepsilon \right)
\end{cases}
\]
Then,
\[
\left \Vert f_{\lambda}\right \Vert _{L^{2}\left(  \mathbb{R}\right)  }^{2}=\int
f_{\lambda}(t)^{2}dt\leq2\varepsilon \lambda^{-1/2}.
\]
Thus,
\begin{equation}
\left \Vert f_{\lambda}\right \Vert _{L^{2}\left(  \mathbb{R}\right)  }%
\lesssim \lambda^{-1/4}. \label{e}%
\end{equation}

For any point $\left(  u,v\right)  $ with $u\in \left(  1-\varepsilon
,1+\varepsilon \right)  $ and $v\in \lambda^{-1/4}\left(  1-\varepsilon
,1+\varepsilon \right)  $, we have $|\lambda \left(  ut^{2}+v^{2}t\right)
D-2D|=|S\left(  u,v\lambda^{1/4},\lambda^{1/2}t\right)  -S(1,1,1)|<\delta$
over the interval $t\in \lambda^{-1/2}\left(  1-\varepsilon,1+\varepsilon
\right)  .$ Thus, the imaginary part of $e^{i\lambda \left(  ut^{2}%
+v^{2}t\right)  }=\sin \left[  \lambda \left(  ut^{2}+v^{2}t\right)  D\right]
>\sin \left(  D\right)  $ , a positive constant not depending on $\lambda.$
So,
\begin{align*}
|T_{\lambda}f_{\lambda}\left(  u,v\right)  |  &  =|\int e^{i\lambda \left(
ut^{2}+v^{2}t\right)  }\Phi(u,v,t)f_{\lambda}(t)dt|\\
&  \geq \sin \left(  D\right)  \int_{t\in \lambda^{-1/2}\left(  1-\varepsilon
,1+\varepsilon \right)  }f_{\lambda}(t)dt\\
&  \gtrsim \lambda^{-1/2}.
\end{align*}
  Thus,$\ $
\begin{align*}
\left \Vert T_{\lambda}f_{\lambda}\right \Vert _{L^{2}\left(  \mathbb{R}%
^{2}\right)  }^{2}  &  \geq \left[  \iint_{\{ \left(  u,v\right)  :u\in \left(
1-\varepsilon,1+\varepsilon \right)  \text{ and }v\in \lambda^{-1/4}\left(
1-\varepsilon,1+\varepsilon \right)  \}}|T_{\lambda}f_{\lambda}\left(
u,v\right)  |^{2}dudv\right] \\
&  \gtrsim \left(  \lambda^{-1/2}\right)  ^{2}\cdot \lambda^{-1/4}%
=\lambda^{-5/4},
\end{align*}
and%
\begin{equation}
\left \Vert T_{\lambda}f_{\lambda}\right \Vert _{L^{2}\left(  \mathbb{R}%
^{2}\right)  }\gtrsim \lambda^{-5/8}. \label{f}%
\end{equation}
By (\ref{e}) and (\ref{f}), we have
\[
\left \Vert T_{\lambda}\right \Vert _{L^{2}\left(  \mathbb{R}\right)
\rightarrow L^{2}\left(  \mathbb{R}^{2}\right)  }\geq \frac{\left \Vert
T_{\lambda}f_{\lambda}\right \Vert _{L^{2}\left(  \mathbb{R}^{2}\right)  }%
}{\left \Vert f_{\lambda}\right \Vert _{L^{2}\left(  \mathbb{R}\right)  }%
}\gtrsim \frac{\lambda^{-5/8}}{\lambda^{-1/4}}=\lambda^{-3/8}.
\]
This proves the sharpness of the decay rate in (\ref{d}).

\end{proof}

\end{document}